\documentclass[12pt]{article}

\usepackage{amssymb,amsthm,amsmath}
\usepackage[usenames,dvipsnames]{xcolor}
\usepackage{soul}
\usepackage{mathabx}
\allowdisplaybreaks
\usepackage{color,graphicx,epsfig}
\usepackage{mathrsfs}
\usepackage{enumerate} 
\usepackage{enumitem} 
\usepackage{esvect}

\newcommand{\comment}[1]{}

\definecolor{teal}{RGB}{0,128,128}
\definecolor{darkpurple}{RGB}{128,0,128}

\usepackage[titletoc,toc,title]{appendix}

\newcommand{\tmod}[1]{{\;\rm (mod\; #1)}}

\newtheorem{theorem}{Theorem}[section]
\newtheorem{prop}[theorem]{Proposition}
\newtheorem{lemma}[theorem]{Lemma}

\theoremstyle{definition}
\newtheorem{defn}[theorem]{Definition}
\newtheorem{rem}[theorem]{Remark}
\theoremstyle{definition}

\def \cB {{\cal B}}

\def \cF {{\cal F}}

\def \cT {{\cal T}}
\def \cU {{\cal U}}

\def \Z {\mathbb Z}

\title{The existence of pyramidal Steiner triple systems over abelian groups}
\author{Yanxun Chang\footnotemark[1],  
Tommaso\ Traetta\footnotemark[2],
Junling Zhou\footnotemark[3]
}
\date{\vspace{-5ex}}

\begin{document}
\maketitle

\footnotetext[1]{School of Mathematics and Statistics, Beijing Jiaotong University, Beijing, 100044, People’s Republic of China. E-mail: yxchang@bjtu.edu.cn}
\footnotetext[2]{DICATAM, Universit\`{a} degli Studi di Brescia, Via Branze 43, 25123 Brescia, Italy. E-mail: tommaso.traetta@unibs.it}
\footnotetext[3]{School of Mathematics and Statistics, Beijing Jiaotong University, Beijing, 100044, People’s Republic of China. E-mail: jlzhou@bjtu.edu.cn}

\begin{abstract} 
  A Steiner triple system STS$(v)$ is called $f$-pyramidal if it has an automorphism group  fixing $f$ points and acting sharply transitively on the remaining $v-f$ points. In this paper, we focus on the STSs that are $f$-pyramidal over some abelian group. Their existence has been settled only for the smallest admissible values of $f$, that is, $f=0,1,3$.
  
  In this paper, we complete this result and determine, for every $f>3$, the spectrum of values $(f,v)$ for which there is an $f$-pyramidal STS$(v)$ over an abelian group.
This result is obtained by constructing difference families relative to a suitable partial spread.
\end{abstract}


\section{Introduction}
  A \emph{Steiner triple system} STS$(v)$ of order $v$ is a pair $(V, \cB)$ where $V$ is a set of $v$ \emph{points} and $\cB$ is a set of unordered \emph{triples} (also called \emph{blocks}) such that any two distinct points belong to exactly one triple. It is well known that an STS$(v)$ exists if and only if $v\equiv 1$ or $3 \pmod{6}$. For a general background on STSs we refer the reader to \cite{ColRo99}.
  
Steiner triple systems have been widely studied over the past 170 years, yet there are still several open questions concerning, for example, the existence of STSs with prescribed symmetries. In this paper, we focus our attention on STSs having an automorphism group $G$ fixing $f$ points and acting sharply transitively on the remaining $v-f$ points. Such an STS will be called \emph{$f$-pyramidal} (over $G$). 
Since every STS($v$) can be considered $v$-pyramidal under the action of any group, we will speak of a non-trivial $f$-pyramidal STS($v$) whenever $f<v$.
We notice that the blocks containing only fixed points form an STS($f$) (see \cite{BuRiTra17}), that is, a subsystem of order $f$ of the original $f$-pyramidal STS$(v)$.
We recall that an STS$(v)$ having a proper subsystem of order $f>0$
(regardless their symmetries) exists if and only if 
$f,v\equiv 1,3 \pmod{6}$ and $f<\frac{v}{2}$ (see \cite{DoWi73}). 
One therefore obtains the following necessary conditions, given in \cite[Lemma 1.1]{BuRiTra17}.

\begin{lemma}[\cite{BuRiTra17}]
\label{nec}
A necessary condition for the existence of 
$f$-pyramidal STS$(v)$ is that  $f=0$ or $f\equiv 1, 3 \pmod{6}$, and $f=v$ or  $f<\frac{v}{2}$.
\end{lemma}

Pyramidal STSs have been deeply studied for the smallest values of $f$, that is, when $f=0,1,3$. It is worth pointing out that an $f$-pyramidal STS is more commonly called \emph{regular} when $f=0$, and \emph{$1$-rotational} when $f=1$.
It was proved in \cite{P} that a \emph{regular} STS$(v)$ exists for every admissible $v$, whereas the spectrum of values $v$ for which there is a $3$-pyramidal STS$(v)$ was completely determined in \cite{BuRiTra17}. On the other hand, although $1$-rotational STSs have been widely investigated in a series of papers \cite{BoBuRiTra12, Bu01, Mi07, PheRo81}, their existence remains an open problem 
whenever $v\equiv 1 \pmod{24}$ and the following conditions simultaneously hold:
$v=(p^3-p)n+1\equiv1$ (mod 96) with $p$ a prime;
$n\not\equiv0$ (mod 4); the odd part of $v-1$ is square-free and without prime factors $\equiv1 \pmod{6}$ (see \cite{BoBuRiTra12}). 
It is worth pointing out that automorphism groups exhibiting a pyramidal action have recently appeared in a slightly different context: for example, a recent classification \cite{PaPeZy} of the full automorphism groups of the five symmetric $(15,8,4)$-designs 
revealed that one of these designs possesses an automorphism group with a 7-pyramidal action on the point-set.

We recall that an STS whose set of triples can be partitioned into parallel classes (where each parallel class is in its turn a partition of the point-set) is called a \emph{Kirkman triple system} (briefly, KTS). The very few results on $f$-pyramidal KTSs are surveyed in \cite{BoBuGaRiTra21} where, in particular, the $3$-pyramidal approach has proven successful to construct the first infinite families of highly symmetric KTSs whose orders fill a congruence class. In the same paper, it is shown that a group having a $3$-pyramidal action on a KTS has necessarily three involutions (i.e., elements of order $2$) and they must be pairwise conjugate. These groups have then been recently classified in \cite{GaoGar23} (see also \cite{GaoGar25}).

In this paper, we determine sufficient and necessary conditions for the existence of an $f$-pyramidal STS$(v)$ over an abelian group, without any restriction on $f$. 
The following theorem settles the existence problem when $0\leq f\leq 3$.

\begin{theorem}[\cite{Bu01, BuRiTra17, P, PheRo81}]\label{f=0,1,3}
  Let $0\leq f \leq 3$. There exists an $f$-pyramidal STS$(v)$ over an abelian group if and only if 
  \begin{enumerate}
    \item $f=0$ and $v\equiv 1,3\pmod{6}$,
    \item $f=1$ and either $v\equiv 3,9 \pmod{24}$ or 
          $v\equiv 1,19 \pmod{72}$,
    \item $f=3$ and either $v\equiv 7,15 \pmod{24}$ or 
          $v\equiv 3,19 \pmod{48}$.
  \end{enumerate}
\end{theorem}

Here is the main result of this paper which completes the previous theorem by determining, for every $f>3$, the spectrum of positive values $(f,v)$ for which there is a non-trivial $f$-pyramidal STS$(v)$ over an abelian group.


\begin{theorem}\label{main} 
  There exists an $f$-pyramidal STS$(v)$, with $3<f<v$, over some abelian group if and only if $f=2^m-1$ for some $m\geq 3$ and one of the following conditions hold:
  \begin{enumerate}[label=$(\alph*)$]
    \item $v \equiv 2^{m+1} -1 \pmod{2^m3}$,
    \item $m$ is even and
    $v \equiv 2^{m} -1 \pmod{2^m3}$,
    \item $m$ is odd and
    $v \equiv 2^{m} -1 \pmod{2^m9}$.
  \end{enumerate}
\end{theorem}

The above result (proven in Section 4) is obtained after constructing, in Section 3 (see Theorem \ref{main:2}), suitable difference families that generate pyramidal STSs (see Theorem \ref{DF}). All preliminary concepts and results are given in the following section.

\section{Preliminaries}
Let $(G,+)$ be a finite group and let $\Sigma$ be a family of subgroups of $G$. We say that $\Sigma$ is a partial spread (briefly, PS) of $G$ if its elements have pairwise trivial intersection; we speak of a PS of type $\tau=\{n_1^{f_1},\ldots, n_t^{f_t}\}$ (briefly, $\tau$-PS) whenever $\Sigma$ has cardinality $f_1+\cdots+f_t$ and contains exactly $f_i$ groups of order $n_i$, for every $1\leq i\leq t$. 

Given a subgroup $H$ of $G$, we denote by $\sigma(H)$ the set of all non trivial subgroups of $H$ of prime order. Clearly, $\sigma(H)$ is a partial spread of $H$.

Given a triple $T=\{a,b,c\}\subset G$, we denote by $\Delta T=\pm\{a-b,a-c,b-c\}$ the 
\emph{list of differences} of $T$. Given a family $\cT$ of unordered triples, the list of differences of $\cT$ is the multiset union $\Delta \cT=\bigcup_{T\in\cT}\Delta T$.

A $(G,\Sigma,3,1)$-difference family (briefly, DF) is a set $\cT$ of unordered triples of $G$ such that $\Delta \cT = G\setminus \bigcup_{S\in\Sigma} S$. We will refer to 
$\cT$ as a $(G,H,3,1)$-DF or $(G,\tau,3,1)$-DF whenever $\Sigma=\{H\}$ or $\Sigma$ is a $\tau$-PS, respectively. Furthermore, when $H$ is isomorphic to the cyclic group $\Z_h$ of order $h$, we more commonly speak of a
$(G,\Z_v,3,1)$-DF. These types of difference families (\emph{relative} to the subgroup $H$ or \emph{relative} to a partial spread, respectively) were introduced by M. Buratti in \cite{Bu98, Bu01_2}. Note that when $\Sigma=\{\{0\}\}$, we simply write 
$(G,3,1)$-DF. Clearly, the union of two difference families can sometimes form a new difference family (see, for example, Lemma \ref{xy} and Proposition \ref{ml=even0}).

\begin{rem}\label{rem} If $H$ is an elementary abelian $p$-group, say $H\simeq Z_{p}^n$, then
$\sigma(H)$ is a $\{p^f\}$-PS of $H$ with 
$f=\frac{p^n-1}{p-1}$. We also note that if $H$ is a subgroup of $G$, then a $(G,H,3,1)$-DF can be viewed as a $(G,\sigma(H),3,1)$-DF.
\end{rem}

The following characterization for $f$-pyramidal STS$(v)$ over an arbitrary group $G$ was proven in \cite{BuRiTra17}. We recall that an {\it involution} of $G$ is an element $x\in G$ of order $2$; hence, $x=-x$.
\begin{theorem}[{\cite[Theorem 2.1]{BuRiTra17}}]\label{DF}
There exists a non-trivial $f$-pyramidal STS$(v)$ under a group $G$ if and only if the following conditions hold:
\begin{enumerate}
  \item $f=0$ or $f\equiv 1, 3 \pmod{6}$, and $f<\frac{v}{2}$,
  \item $G$ has order $v-f$ and contains exactly $f$ involutions, and
  \item there exists a $(G,\{2^f,3^e\},3,1)$-DF for some $e\geq 0$.
\end{enumerate}
\end{theorem}

\begin{rem} It is not difficult to check that given a $(G,\Sigma,3,1)$-DF relative to
a partial spread $\Sigma$ of type $\{2^f, 3^e\}$, we have that $\Sigma$ must contain all subgroups of $G$ of order $2$, that is, $f$ must coincide with the number of involutions of $G$. Therefore, in condition (2) of the previous theorem it is enough to require that $|G|=v-f$.
\end{rem}

In this paper, we focus our attention on $f$-pyramidal STSs over abelian groups.
The following result provides some necessary conditions for the existence of a DF as in 
Theorem \ref{DF}.(3).

\begin{lemma}\label{lemma_1} 
Let $G$ be an abelian group having exactly $f>0$ involutions. If there exists a  $(G,\{2^f,3^e\},3,1)$-DF
for some $e\geq 0$, then the following conditions hold:
\begin{enumerate}
\item $f=2^m-1$ for some integer $m\geq 1$;
\item $|G|=2^{m+\ell}d$ for some odd integer $d\geq 1$ and an integer $\ell\geq 0$; 
\item $G=\Z_{2^{\alpha_1}}\times \cdots \times \Z_{2^{\alpha_m}} \times H$ for a suitable abelian group $H$ of order $d$
and for suitable positive integers $\alpha_1, \ldots, \alpha_m$ such that 
$\alpha_1+\cdots + \alpha_m=m+\ell$;
\item If $3\mid d$, then 
$e\equiv (-1)^m \pmod{3}$;
  \item $d\equiv (-1)^\ell \pmod{6}$
  or  $d\equiv
 \begin{cases}
   3 \pmod{6} & \text{if $m$ is even},\\
   9 \pmod{18} & \text{if $m$ is odd}.
 \end{cases}$ 
\end{enumerate}
\end{lemma}
\begin{proof} 
Let $G$ be an abelian group with exactly $f>0$ involutions 
 and order $2^nd$, where $n>0$ and $d\geq 1$ is odd. 
Then, $G=P\times H$ where $P$ is the Sylow 2-subgroup of $G$ and $H$ has order $d$.
By the fundamental theorem of finite abelian groups, it follows that
$P=\bigtimes_{i=1}^{m} \Z_{2^{\alpha_i}}$
where $m\geq 1$, each $\alpha_i\geq 1$ and $\sum_{i=1}^m \alpha_i=n$. 

Since $G$ is abelian, the set $I(G)$ containing all involutions of $G$ and the zero element coincides with the unique elementary abelian 2-subgroup of $P$, that is,
\[ I(G) = 2^{\alpha_1-1}\Z_{2^{\alpha_1}}\times \cdots \times 2^{\alpha_m-1}\Z_{2^{\alpha_m}}.
\]
Hence, $f=|I(G)|-1 = 2^m-1$. We have therefore proven items (1), (2) and (3), with
$\ell=n-m$.

Now, let $\cF$ be a $(G,\Sigma,3,1)$-DF where $\Sigma$ is a $\{2^f,3^e\}$-PS.
Since $|\Delta \cF| =6 |\cF|=|G\setminus\bigcup_{S\in \Sigma} S|$ and 
$|\bigcup_{S\in \Sigma} S| = 2^m+2e$, we have that
  \[|G\setminus\textstyle{\bigcup_{S\in \Sigma}} S|=2^{m+\ell}d -(2^m+2e) = 2^{m}(2^{\ell}d-1)-2e\equiv 0 \pmod{6},\]
  hence $2^{m-1}(2^{\ell}d-1) \equiv e \pmod{3}$.
  If we assume that $3|d$, then $2^{\ell}d \equiv 0\pmod{3}$, hence $e \equiv 2^{m} \pmod{3}$, thus proving item $(4)$. If, in addition, $m$ is odd, then $e\geq 2$, that is, 
  $\Sigma$ contains at least two distinct groups of order $3$. Hence $\Z_3\times \Z_3$ is a subgroup of $G$ (since $G$ is abelian), therefore $9$ divides $|G|$. In other words,
\begin{equation}\label{a}
  \text{if $3|d$ and $m$ is odd, then $9|d$}.
\end{equation}  
  By Theorem \ref{DF}, a $(G,\{2^f,3^e\},3,1)$-DF determines an $f$-pyramidal STS of order $|G|+f = 2^{m+\ell}d+2^m-1 = 2^m(2^\ell d+1) -1 \equiv 1,3 \pmod{6}$, hence $2^m(2^\ell d+1) \not\equiv 0 \pmod{6}$. Then, $2^{m-1}(2^\ell d+1) \not\equiv 0\pmod{3}$, that is, 
\begin{equation*}
  2^\ell d\equiv 0,1 \pmod{3}.
\end{equation*}
By taking into account that $d$ is odd, one can easily check that
\begin{enumerate}
\item[$a$.] 
$2^\ell d\equiv 1 \pmod{3} \;\Leftrightarrow\; d\equiv (-1)^\ell \pmod{6}$;
\item[$b$.] 
$2^\ell d\equiv 0 \pmod{3}$ 
$\;\overset{\eqref{a}}{\Leftrightarrow}\; 
d\equiv
 \begin{cases}
   3 \pmod{6} & \text{if $m$ is even},\\
   9 \pmod{18} & \text{if $m$ is odd},
 \end{cases}$
\end{enumerate}
thus proving item $(5)$.
\end{proof}

At this stage, a reader exclusively interested in the \textit{only if} part of Theorem \ref{main} can skip the next few pages and go directly to Section \ref{section:main}.\\

In the following, we recall some results on abelian difference families, and difference matrices that we will need in our constructions.

\begin{lemma}[\cite{P}]\label{cyclic_DF}
  There exists a $(\Z_{v},3,1)$-DF for every $v\equiv 1 \pmod{6}$. 
  There exists a $(\Z_{v},\Z_3, 3,1)$-DF for every $v\equiv 3 \pmod{6}$, with $v\neq 9$,
  and a $(\Z_{3}^2,\Z_3, 3,1)$-DF.
\end{lemma}

\begin{lemma}[{\cite[Theorem 2.2]{PheRo81}}]\label{1rot_DF}
  There is a  $(\Z_{v}, \Z_2, 3,1)$-DF if and only if
  $v\equiv 2,8 \pmod{24}$.
\end{lemma}

The following two results follow from 
\cite[Theorem 5.1]{Bu01} and \cite[Theorem 4.4]{BuRiTra17}, respectively.
\begin{lemma}[\cite{Bu01}]\label{1rot_DF_abelian}
  There is a  $(\Z_{3}\times \Z_{6n}, \{2^1, 3^2\}, 3,1)$-DF whenever $n\equiv 0,1 \pmod{4}$.
\end{lemma}

\begin{lemma}[{\cite{BuRiTra17}}]\label{3pyr_DF}\hspace*{\fill} 
  \begin{enumerate}
    \item For every $n\geq 1$, there exists a  $(\Z_{4}\times\Z_{12n}, \Sigma, 3,1)$-DF where $\Sigma$ is a partial spread of type $\{2^3,3\}$.

    \item For every $n\geq 0$, there exists a  
 $(\Z_{2}^2 \times \Z_{6n+3}, \Sigma, 3,1)$-DF where $\Sigma$ is a partial spread of type $\{2^3,3\}$.
 \end{enumerate}
\end{lemma}

We recall that a $(G,3,1)$-difference matrix (DM) is a $3\times |G|$ array, with entries from an additive group $G$, such that the difference of any two distinct rows is a permutation of $G$.

\begin{lemma}[\cite{Hall-Paige}, {\cite[Theorem 17.9]{DM}}]\label{DM}
  There exists a $(G,3,1)$-DM if and only if the Sylow $2$-subgroups of $G$ are trivial or noncyclic.
\end{lemma}

The following is a natural generalization of a standard method of expanding difference families by using difference matrices. 
\begin{lemma}
\label{composition}
  If there exists an $(H,\Sigma,3,1)$-DF and a $(K,3,1)$-DM, then there exists a $(H\times K, \Sigma',3,1)$-DF, where $\Sigma'=\{S\times K\mid S\in \Sigma\}$.
\end{lemma}
\begin{proof}
  Let $\cT$ be an $(H,\Sigma,3,1)$-DF and let $M=(m_{ij})$ be a $(K,3,1)$-DM. For every $T=\{a,b,c\}\in \cT$, we define the triple $T_j=\{(a,m_{1j}), (b,m_{2j}), (c,m_{3j})\}$
  and let $\cT'=\{T_j\mid T\in \cT, 1\leq j\leq |K|\}$. We notice that
\[\Delta T_j= \pm \{(b-a, m_{2j}-m_{1j}), (c-b, m_{3j}-m_{2j}), (a-c, m_{1j}-m_{3j})\}.\]
By definition of a difference matrix, it follows that
\[\bigcup_{j=1}^{|K|}\Delta T_j= 
(\pm \{b-a, c-b, a-c\})\times K = \Delta T \times K.\]
Therefore, 
\[\Delta \cT' = \bigcup_{T\in \cT} (\Delta T \times K) = (\bigcup_{T\in \cT} \Delta T) \times K = (H\setminus \bigcup_{S\in\Sigma} S) \times K = 
(H\times K)\setminus \bigcup_{S\in\Sigma} (S \times K).\]
Hence, $\cT'$ is the desired DF.
\end{proof}

In our constructions, we will make use of Langford sequences \cite{Langford} whose definition is recalled in the following.

\begin{defn}
Let $k, a$ and $b$ be non-negative integers with $1\leq k \leq 2a+1$ and $b\geq 1$.
A {\it $k$-extended Langford sequence} of {\it order $a$} and {\it defect $b$} is a sequence of $a$ integers 
$(s_1, s_2, \ldots, s_a)$ such that
\[\{s_i, s_i +i + (b-1) \;|\; i=1, \ldots, a\} = [1,  2a+1]\setminus \{k\}.\] 
If $b=1$, one speaks of a {\it $k$-extended Skolem sequence}.
\end{defn}

\begin{defn} Let $k, a$ and $b$ be non-negative integers with $1\leq k \leq 2a+1$ and $b\geq 1$.
  We say that a triple $(k,a,b)$ is {\it Langford admissible} if either $(a,k)=(0,1)$ or $2b-1 \leq a$, and one of the following conditions holds:
  \begin{enumerate}
  \item $b$ is odd,  $k$ is odd,  $a \equiv 0,1 \tmod{4}$; 
  \item $b$ is odd,  $k$ is even, $a \equiv 2,3 \tmod{4}$; 
  \item $b$ is even, $k$ is odd,  $a \equiv 0,3 \tmod{4}$;  
  \item $b$ is even, $k$ is even, $a \equiv 1,2 \tmod{4}$.
\end{enumerate}
\end{defn}

The following is a restricted version of a more general result on the existence of Skolem and Langford sequences.

\begin{theorem}[\cite{Baker95, LiJi98, LiMor04}]
\label{lang} 
If $b\in[1,4]$ and $(k,a,b)$ is Langford admissible, then  there exists a $k$-extended Langford sequence of order $a$ and defect $b$.
\end{theorem}

The following is a standard method, based on Langford sequences, to build a family of triples of $\Z$ with a given list of differences.

\begin{lemma} \label{lang_triples}
If $(k,a,b)$ is Langford admissible and $b\in[1,4]$, then
there exists a set $\cT$ of triples of $\Z$, with $|\cT|=a$, such that
\[\Delta \cT = \pm ([b, 3a+b]\setminus \{k+a+b-1\}).
\]
\end{lemma}
\begin{proof} By Theorem \ref{lang}, there exists a $k$-extended Langford sequence of order $a$ and defect $b$, say $(s_1, \ldots, s_a)$. 
It is then enough to take $\cT=\{T_i\mid 1\leq i\leq a\}$, where
$T_i=\{0, s_i + a+b-1, s_i + i + a+ 2(b-1)\}$, for $1\leq i\leq a$.
\end{proof}

\section{The existence of a $(G,\{2^f,3^e\},3,1)$-DF over an abelian group $G$}
In this section, we show that the converse of Lemma \ref{lemma_1} holds for at least one abelian group of each admissible order, provided that $f\geq 7$. In other words, we show the following result whose proof is split into Sections \ref{sec:3.1} and \ref{sec:3.2}.
\begin{theorem}\label{main:2} Let $d,\ell,m$ be integers with
$d\geq 1$ odd, $\ell\geq 0$ and $m\geq 3$. If 
\[
\text{$d\equiv (-1)^\ell \pmod{6}$\;\;\; or\;\;\;$
d\equiv
 \begin{cases}
   3 \pmod{6} & \text{if $m$ is even},\\
   9 \pmod{18} & \text{if $m$ is odd},
 \end{cases}$}
\]
then there exists an abelian group $G$ of order $2^{\ell+m} d$, satisfying the following conditions:
  \begin{enumerate}
  \item $G$ has exactly $f=2^{m}-1$ involutions, and
  \item there is a $(G,\{2^f,3^e\},3,1)$-DF, for some $e\in \{0,1,2\}$.
\end{enumerate}
\end{theorem}

\subsection{Case $d\equiv (-1)^\ell \pmod{6}$}
\label{sec:3.1}
In this subsection, we prove Theorem \ref{main:2} under the assumption that $d\equiv (-1)^\ell \pmod{6}$; hence, we build a $(G,\{2^f\},3,1)$-DF over a suitable abelian group $G$ of order $2^{\ell+m}d$ having $f=2^m-1$ involutions. In particular, the cases $\ell=0$ and $\ell\geq 2$ are dealt with in Propositions \ref{1:l=0} and \ref{1:l>=2}, respectively, whereas the case $\ell=1$ is solved in Propositions \ref{1:l=1,m<>4} and \ref{1:l=1,m=4}.

\begin{prop}\label{1:l=0}
  There is a $(\Z_2^m \times \Z_d,\{2^f\},3,1)$-DF whenever $d\equiv 1 \pmod{6}$ and $m\geq 2$.
\end{prop}
\begin{proof} Set $G=\Z_2^m \times \Z_d$, where $d \equiv 1 \pmod{6}$. By Lemma \ref{cyclic_DF}, there exists a $(\Z_{d},3,1)$-DF.
Also, by Lemma \ref{DM}, there exists a $(\Z_2^m, 3, 1)$-DM.
Therefore, Lemma \ref{composition} guarantees the existence of a
 $(G, \Z_2^m \times\{0\}, 3,1)$-DF, or equivalently,
(see Remark \ref{rem}) a $(G, \sigma(H), 3,1)$-DF, 
where $\sigma(H)$ is a $\{2^{f}, 3^0\}$-PS, with $f=2^{m}-1$.
\end{proof}

\begin{prop}\label{1:l>=2}
Let $G=\Z_2^{m-1} \times \Z_{2^{\ell+1}d}$. 
There is a $(G,\{2^f\},3,1)$-DF whenever $\ell\geq 2, m\geq 3$ and $d\equiv (-1)^\ell \pmod{6}$.
\end{prop}
\begin{proof}
  Set $H'=2^\ell d\Z_{2^{\ell+1}d}\simeq \Z_2$ and 
  $H=\Z_2^{m-1}\times H'\simeq \Z_2^{m}$.
  
Since $\ell\geq 2$ and $d\equiv (-1)^\ell \pmod{6}$, we have that $2^{\ell+1}d\equiv 8 \pmod{24}$. Therefore, by Lemma \ref{1rot_DF}, there is a $(\Z_{2^{\ell+1}d}, H', 3,1)$-DF. Also, by Lemma \ref{DM}, there exists a $(\Z_2^{m-1}, 3, 1)$-DM.
Then, Lemma \ref{composition} produces
a $(G, H, 3,1)$-DF which is equivalent 
(see Remark \ref{rem}) to a $(G, \sigma(H), 3,1)$-DF, 
where $\sigma(H)$ is $\{2^{f}\}$-PS.
\end{proof}

\begin{lemma}\label{mld=315} 
There exists a 
$(\Z_2^{2}\times \Z_{4d}, H,3,1)$-DF, with $H\simeq \Z^3_{2}$,
whenever $d\equiv 5\pmod{6}$.
\end{lemma}
\begin{proof} Set $\Z_2^2 = \{0, \alpha,\beta,\gamma\}$ so that $\alpha,\beta,\gamma$ are involutions and $\alpha+\beta=\gamma$.

We start by defining a set $\cT=\{T_1, \ldots, T_{2d-3}\}$ of 
$2d-3$ triples as follows:

\[T_i = 
  \begin{cases}
    \{(0, 0), (\alpha,i), (\gamma,2i)\} & 
      \text{if $1\leq i\leq d-1$},\\
    \{(0, 0), (\alpha,i+1), (\gamma,2i+1)\} & 
      \text{if $d\leq i\leq 2d-3$}.
  \end{cases}      
\]   
Also, let $\cT'=\{T'_\alpha, T'_\beta, T'_\gamma\}$ where:
\begin{align*}
 T'_\alpha &= \{(0, 0), (0, d+1), (\alpha, 2d+1)\},\\
 T'_\beta  &= \{(0, 0), (0, 1), (\beta, 2d+2)\},\\
 T'_\gamma &= \{(0, 0), (0, 2), (\gamma, 3)\}.
\end{align*}
Note that 
$\Delta \cT = \bigcup_{x\in\{\alpha,\beta,\gamma\}} \{x\}\times D_x$ and $\Delta \cT' = \bigcup_{x\in\Z_2^2} \{x\}\times D'_x$, 
where
\begin{align*}
  D_\alpha &=\pm[1,d-1]\,\cup\,\pm [d+1, 2d-2] = \Z_{4d}\setminus(\pm\{0,d, 2d-1, 2d\}),\\
  D_\beta &=\pm[1,2d-3]=\Z_{4d}\setminus(\pm\{0,2d-2, 2d-1, 2d\}),\\  
  D_\gamma &=\{\pm2\}\,\cup\, \pm[4,2d-1] = \Z_{4d}\setminus(\pm\{0,1, 3, 2d\}),\\
  D'_0  &=\pm\{1,2,d+1\},\\
  D'_\alpha &=\pm\{d, 2d+1\}=\pm\{d, 2d-1\},\\
  D'_\beta &=\pm\{2d+1, 2d+2\} = \pm\{2d-2, 2d-1\},\\
  D'_\gamma &=\pm\{1, 3\}.  
\end{align*}
Hence, 
$\Delta (\cT\,\cup\,\cT') = \pm \big(\{0\}\times \{1,2,d+1\}\big)
\,\cup\, 
\big(\{\alpha,\beta,\gamma\}\times \Z_{4d}\setminus\{0,2d\}\big)
$. Therefore, it is left to construct a set $\cT''$ of $\frac{2d-4}{3}$ triples of 
$\{0\}\times \Z_{4d}$ such that 
\[\Delta \cT'' = \pm \big(\{0\}\times ([3, 2d-1]\setminus\{d+1\})\big).\]
If $d=5$, we take $\cT''=\{\{(0,0),(0,3),(0,8)\}, \{(0,0),(0,4),(0,11)\}\}$. Otherwise, for $d\geq 11$, set $(k,a,b)=(\frac{d+1}{3},\frac{2d-4}{3},3)$. Since $(k,a,b)$ is Langford admissible,  Lemma \ref{lang_triples} guarantees the existence of the desired set of triples $\cT''$. Therefore, $\cT\,\cup\,\cT'\,\cup\,\cT''$ is
a $(\Z_2^{2}\times \Z_{4d}, H,3,1)$-DF where $H=\Z_2^{2}\times 2d\Z_{4d}\simeq \Z_2^3$.
\end{proof}

\begin{prop}\label{1:l=1,m<>4}
There exists a 
$(\Z_2^{m-1}\times \Z_{4d}, H,3,1)$-DF, with $H\simeq\Z_2^{m}$, whenever  $3\leq m\neq 4$ and $d\equiv 5\pmod{6}$.
\end{prop}
\begin{proof}
By Lemma \ref{mld=315}, there exists a $(\Z_2^{2}\times \Z_{4d}, H',3,1)$-DF, with 
$H'\simeq \Z_2^{3}$. Also, by Lemma \ref{DM}, there exists
a $(\Z_2^{m-3},3,1)$-DM. Then, Lemma \ref{composition} guarantees the existence of
$(\Z_2^{m-1}\times \Z_{4d}, H,3,1)$-DF, where $H=\Z_2^{m-3} \times H'\simeq \Z_2^{m}$.
\end{proof}

We now deal with the missing case in Proposition \ref{1:l=1,m<>4}, that is, $m=4$.
\begin{prop}\label{1:l=1,m=4}
There exists a 
$(\Z_2^{3}\times \Z_{4d}, \Z_2^{4},3,1)$-DF whenever $d\equiv 5\pmod{6}$.
\end{prop}
\begin{proof} 
  Set $\Z_2^2 = \{0, \alpha,\beta,\gamma\}$ so that $\alpha,\beta,\gamma$ are involutions and $\alpha+\beta=\gamma$, and
  set
  $G=\Z_2^2 \times \Z_2 \times \Z_{4d}$.

  Let $\cT=\{T_{i,j} \mid i\in[1,4], j\in[1,d-1]\}$ be the set of triples of $G$ defined below:
  \begin{align*}
    T_{1,j}    &=\{(0,0,0), (\alpha,0,j), (\gamma, 0, 2j+1)\},\\
    T_{2,j}    &=\{(0,0,0), (\alpha,0,j+d), (\gamma, 1, 2(j+d))\},\\
    T_{3,j}    &=\{(0,0,0), (\alpha,1,j), (\gamma, 0, 2j)\},\\
    T_{4,j}    &=\{(0,0,0), (\alpha,1,j+d-1), (\gamma, 1, 2(j+d)-1)\}.                        
  \end{align*}
For every $i\in[1,4]$ and $j\in[1,d-1]$, we have that
\[
\Delta T_{i,j}=\pm
\begin{cases}
  \{(\alpha,0,j),(\beta,0,j+1),(\gamma,0,2j+1)\} & \text{if $i=1$},\\
  \{(\alpha,0,j+d),(\beta,1,j+d),(\gamma,1,2(j+d))\} & \text{if $i=2$},\\
  \{(\alpha,1,j),(\beta,1,j),(\gamma,0,2j)\} & \text{if $i=3$},\\
  \{(\alpha,1,j+d-1),(\beta,0,j+d),(\gamma,1,2(j+d)-1)\} & \text{if $i=4$}.\\      
\end{cases}
\]
Therefore,
\begin{equation}\label{eq1:mld=415(mod12)_first}
\begin{aligned}
\Delta \cT = \bigcup_{x\in\{\alpha, \beta, \gamma\}} 
\{x\}\times \big((\Z_2\times \Z_{4d})\setminus (D_x\cup H)\big)\,\,\,\text{where}\\
\begin{aligned}
  D_\alpha = &\pm \{(0,d), (1,2d+1)\},\\
  D_\beta  = &\pm\{(0,1), (1,d)\},\\
  D_\gamma = &\pm\{(0,1),(1,1)\}, \text{and}\\
  H= &\; \Z_2 \times 2d\Z_{4d}=\{(0,0), (0,2d),(1,0), (1,2d)\}.
\end{aligned}
\end{aligned}
\end{equation}

Now, let $\cT'=\{T'_1, T'_2, T'_3\}$, where
  \begin{align*}
    T'_1   &=\{(0,0,0), (0,1,d+1), (\alpha, 1, 2d+1)\},\\
    T'_2   &=\{(0,0,0), (0,1,d-1), (\beta, 1, d)\},\\
    T'_3   &=\{(0,0,0), (0,1,2), (\gamma, 1, 1)\}.        
  \end{align*}  
Note that 
\begin{equation}\label{eq1:mld=415(mod12)_second}
\begin{aligned}
\begin{aligned}
\Delta\cT' = (\{0\}\times D_0)\,\cup\, \
\bigcup_{x\in\{\alpha, \beta, \gamma\}} \{x\}\times D_x,\;\;\; \text{where}\\
D_0=\pm\{(1,2),(1,d-1),(1,d+1)\}.
\end{aligned}
\end{aligned}
\end{equation}
It is then left to construct a set $\cT''$ of triples such that
\begin{equation}\label{eq1:mld=415(mod12)}
\begin{aligned}
  \Delta \cT'' &=\{0\}\times \big((\Z_2\times \Z_{4d})\setminus (D_0\cup H)\big).
\end{aligned}
\end{equation}
Indeed, letting $\cT^*= \cT
\,\cup\,\cT'\,\cup\,\cT''$, by \eqref{eq1:mld=415(mod12)_first},
\eqref{eq1:mld=415(mod12)_second}
and
\eqref{eq1:mld=415(mod12)}, it follows that 
$\Delta \cT^* = G\setminus (\Z_2^2\times H)$. Since $\Z_2^2\times H \simeq \Z_2^4$, we have that $\cT^*$ is the desired DF.

We first consider the set $\cU=\{U_i\mid 1\leq i\leq d\}$ of triples defined as follows:
\begin{align*} 
     U_1   &=
     \begin{cases}
       \{(0,0,0), (0,1,d+3), (0, 0, -1)\} &\text{if $d\equiv 5 \pmod{12}$},\\
       \{(0,0,0), (0,1,d+2), (0, 0, -2)\}, & \text{if $d\equiv 11 \pmod{12}$},   
     \end{cases}\\      
     U_2   &= 
     \begin{cases}
       \{(0,0,0), (0,1,1), (0, 0, -2)\}, & \text{if $d\equiv 5 \pmod{12}$},\\
       \{(0,0,0), (0,1,1), (0, 0, 4)\}, & \text{if $d\equiv 11 \pmod{12}$},       
     \end{cases}\\
     U_3   &=
     \begin{cases}
       \{(0,0,0), (0,1,d), (0, 0, 2d+2)\}, & \text{if $d\equiv 5 \pmod{12}$},\\
       \{(0,0,0), (0,1,d), (0, 0, -3)\}, & \text{if $d\equiv 11 \pmod{12}$},\\       
     \end{cases} \\              
     U_i   &=\{(0,0,0), (0,1,i), (0, 1, 2d+3-i)\},\, \text{for $4\leq i\leq d-2$,}\\
     U_{d-1} &=
     \begin{cases}
       \{(0,0,0), (0,0,3), (0, 0, 2d-3)\}, & \text{if $d\equiv 5 \pmod{12}$},\\
       \{(0,0,0), (0,0,2d-6), (0, 0, 2d-1)\}, & \text{if $d\equiv 11 \pmod{12}$},\\        
     \end{cases}\\     
     U_{d}   &=
     \begin{cases}
       \{(0,0,0), (0,0,5), (0, 0, 2d+1)\}, & \text{if $d\equiv 5 \pmod{12}$},\\
       \{(0,0,0), (0,0,2d-3), (0, 0, 2d-2)\}, & \text{if $d\equiv 11 \pmod{12}$}.          
     \end{cases}         
\end{align*}
Note that $\Delta \cU = \pm\big(\{(0,0)\}\times E\big)\,\cup\, 
\pm\big(\{(0,1)\}\times ([1,2d-1]\setminus\{2,d-1,d+1\})\big)$,
where 
\begin{align*}
E =\, &\{2j+1\mid j\in [0,d-1]\}\,\cup\\
&\begin{cases}
  \{2,2d-6,2d-4,2d-2\} & \text{if $d\equiv 5 \pmod{12}$},\\
  \{2,4,2d-6,2d-2\} & \text{if $d\equiv 11 \pmod{12}$}.
\end{cases}
\end{align*}
To construct the remaining triples we use Langford sequences.
\begin{enumerate}
\item
If $d\equiv 5\pmod{12}$, the triple $(k,a,b) = (\frac{2}{3}(d-5)+1,\frac{d-5}{3},2)$,
is Langford admissible. Therefore, Lemma \ref{lang_triples} guarantees the existence
of a set $\cU'$ of triples of $\{(0,0)\}\times \Z$ such that
$\Delta \cU' = \{(0,0)\}\times [2, d-4]$. Note that ${\cal U}'=\emptyset$ when $d=5$.

\item
If $d\equiv 11\pmod{12}$, the triple $(k,a,b) = (\frac{2}{3}(d-5),\frac{d-5}{3},3)$,
is Langford admissible, when $d> 11$.
Therefore, Lemma \ref{lang_triples} guarantees the existence
of a set $\cU'$ of triples of $\{(0,0)\}\times \Z$ such that
$\Delta \cU' = \{(0,0)\}\times ([3, d-2]\setminus\{d-3\})$.
\end{enumerate}
In both cases, we have that $\cT'' = \cU\,\cup\, 2\cdot\cU'$ satisfies \eqref{eq1:mld=415(mod12)}.\\

It is left to deal with the case $d=11$. Let $\cT''=\{T_1, \ldots, T_{13}\}$ be the set of triples defined below:
\begin{align*}
  T_1 &= \{(0,0,0), (0,0,1), (0, 0, 21)\},\\
  T_2 &= \{(0,0,0), (0,0,3), (0, 0, 19)\},\\
  T_3 &= \{(0,0,0), (0,0,2), (0, 0, 8)\},\\
  T_4 &= \{(0,0,0), (0,0,4), (0, 0, 18)\},\\
  T_5 &= \{(0,0,0), (0,1,1), (0, 0, 12)\},\\
  T_6 &= \{(0,0,0), (0,1,23), (0, 0, 10)\},\\
  T_i &= \{(0,0,0), (0,1,i-4), (0, 1, 27-i)\},\, \text{for $7\leq i\leq 13$.}
\end{align*}
One can check that $\cT''$ satisfies  \eqref{eq1:mld=415(mod12)}.
\end{proof}

\subsection{Cases $d\equiv 3\pmod{6}$ 
and $d\equiv 9 \pmod{18}$ 
}
\label{sec:3.2}
In this subsection, we build a $(G,\{2^f, 3^e\},3,1)$-DF, with $e \geq 1$, over a suitable abelian group $G$ of order $2^{\ell+m}d$ having $f=2^m-1$ involutions, under the assumption that
\[\text{$d\equiv 9 \pmod{18}$ when $m$ is odd, otherwise $d\equiv 3 \pmod{6}$.}\]
In particular, the cases $\ell=0, 1$  and $\ell\geq 2$ are dealt with in 
Propositions \ref{ml=even0}, \ref{ml=even1}, \ref{ml=even>1}, when $m\geq4$ is even, and in Propositions \ref{ml=odd0}, \ref{ml=odd1}, \ref{ml=odd>1}, when $m\geq3$ is odd.

\begin{lemma} \label{mld=413}
  There exists a $(\Z^{3}_{2} \times \Z_{12},\{2^{15}, 3\},3,1)$-DF.
\end{lemma}
\begin{proof}  
  Set $\Z_2^2 = \{0, \alpha_1,\alpha_2,\alpha_3\}$ so that $\alpha_1,\alpha_2,\alpha_3$ are involutions and $\alpha_1+\alpha_2=\alpha_3$.
   Also, set
  $G=\Z_2^2 \times H$, where $H=\Z_2 \times\Z_{12}$, and  $K=\Z_2 \times 6\Z_{12}=\{(0,0),(0,6),(1,0),(1,6)\} \simeq \Z_2^2$. 
  
  Let $\cT=\{T_1,\ldots, T_{13}\}$ be the set of triples defined below, where $\alpha_4=\alpha_1$:
  \begin{align*}
    T_i    &=\{(0,0,0), (\alpha_i,1,4), (\alpha_{i+1}, 1, 9)\},\; \text{for $1 \leq i \leq 3$,}\\
    T_{i+3}&=\{(0,0,0), (\alpha_i,1,2), (\alpha_{i+1}, 1, 5)\},\, \text{for $1 \leq i \leq 3$,}\\
    T_7    &=\{(0,0,0), (\alpha_1,0,4), (\alpha_3,0,8)\},\\
    T_8    &=\{(0,0,0), (\alpha_1,0,1), (\alpha_3,0,2)\}\},\\
    T_9    &=\{(0,0,0), (0,1,1), (\alpha_1,1,11)\},\\
    T_{10} &=\{(0,0,0), (0,1,3), (\alpha_2,1,1)\},\\
    T_{11} &=\{(0,0,0), (0,1,2), (\alpha_3,1,1)\},\\
    T_{12} &=\{(0,0,0), (0,1,4), (0,0,11)\},\\
    T_{13} &=\{(0,0,0), (0,0,3), (0,0,10)\}.              
  \end{align*}
One can check that
$\bigcup_{i=1}^{11}\Delta T_i = \big(\{0\}\times D\big)\,\cup\, \big(\{\alpha_1,\alpha_2,\alpha_3\}\times (H\setminus K)\big)$, and
$\Delta T_{12}\,\cup\, \Delta T_{13}=\{0\}\times D'$, where
\begin{align*}
  D =&\pm \{(1,1), (1,2), (1,3)\},\\
  D'=&\pm\{(0,1), (0,2), (0,3), (0,5), (1,4), (1,5)\}.
\end{align*}
Notice that $D\,\cup\, D' = H\setminus (K\, \cup\, \{\pm (0,4)\})$.
Therefore, $\Delta\cT = G\setminus S$, where
$S= (\Z_2^2\times K)\,\cup\, \{\pm(0,0,4)\}$, 
hence $\cT$ is the desired difference family.
\end{proof}

\begin{lemma} \label{xy}
Let $x\geq 2y$ be an even integer, with $y=1,2$. 
Then,
  there exists a $(\Z^{x-1}_{2} \times \Z_{2^y}\times \Z_{3}, \{2^{2^x-1},3\},3,1)$-DF.
\end{lemma}
\begin{proof}
We proceed by induction on the even values of $x\geq 2y$. If $x=2y$, the result follows from Lemma \ref{3pyr_DF} when $y=1$ and Lemma \ref{mld=413} when $y=2$. Assume the assertion holds for $x\geq 2y$ even. This is equivalent to saying that 
there is a $(G, \{H',K'\},3,1)$-DF, say $\cT$, where $G=\Z^{x-1}_{2} \times \Z_{2^y}\times \Z_{3}$,
\begin{align*}
  H' =&\; \Z^{x-1}_{2} \times 2^{y-1}\Z_{2^y} \times \{0\}\simeq \Z^{x}_{2},\; \text{and}\;
  K' = \{(0,0)\}\times \Z_3\simeq \Z_3.
\end{align*} 
Indeed, $\Sigma'=\sigma(H')\,\cup\,\{K'\}$ is a $\{2^{2^x-1},3\}$-PS, and $\cT$ can be seen as a DF relative to $\Sigma'$.

Since, by Lemma \ref{DM}, there is $(\Z^{2}_{2},3,1)$-DM, we use 
Lemma \ref{composition} to obtain a
$(\Z^{2}_{2} \times G, \{H,K\},3,1)$-DF, say $\cU$,
where 
\[ H=\Z_{2}^2\times H'\simeq \Z^{x+2}_{2}
,\; \text{and}\;
  K = \Z_{2}^2 \times K'\simeq \Z_{2}^2 \times \Z_3.
\]
Finally, by Lemma \ref{3pyr_DF}, there is a 
\[\big(K, \big\{\Z_{2}^2 \times\{0\}, \{0\}\times K'\big\},3,1\big)\text{-DF},\] 
say $\cU'$.
Therefore, $\cU\,\cup\,\cU'$ is a 
$(\Z^{2}_{2} \times G, L,3,1)$-DF, where 
\[L=\{H, \Z_{2}^2 \times\{0\},\{0\}\times K'\},\] 
or equivalently,
a $(\Z^{2}_{2} \times G, \Sigma,3,1)$-DF, where 
\[\Sigma = \sigma(H)\,\cup\,\sigma(\Z_{2}^2 \times\{0\})\,\cup\,\{\{0\}\times K'\}\] 
is a $\{2^{2^{x+2}-1}, 3\}$-PS. Since $\Z^2_2\times G \cong \Z^{x+1}_{2} \times \Z_{2^y}\times \Z_{3}$, we have proven the induction step.
\end{proof}

\begin{prop}\label{ml=even0}
There exists a $(\Z^{m}_{2}\times \Z_{d},\{2^{2^m-1}, 3\},3,1)$-DF for every even $m\geq 4$ and $d\equiv 3 \pmod{6}$.
\end{prop}
\begin{proof}
Set $G=\Z^{m-2}_{2} \times H$, where $H=\Z_2^2 \times \Z_{d}$.
By Lemma \ref{3pyr_DF}, there exists an
$(H,\{K_1, K_2\},3,1)$-DF where $K_1=\Z^{2}_{2}\times \{0\}$ and $K_2=\{0\} \times \frac{d}{3}\Z_d\simeq \Z_3$. 
By Lemma \ref{DM}, there exists a $(\Z_2^{m-2}, 3, 1)$-DM.
Then, by Lemma \ref{composition}, we get a $(G, \{G_1, G_2\}, 3,1)$-DF, say $\cT$, 
where 
\begin{align*}
  G_1 &=\Z_2^{m-2}\times K_1 = \Z^{m}_{2}\times\{0\},\;\text{and} \\
  G_2 &=\Z_2^{m-2}\times K_2 = \Z^{m-2}_{2} \times\{0\}\times \textstyle{\frac{d}{3}}\Z_d \simeq \Z^{m-2}_{2} \times \Z_3.
\end{align*}
By Lemma \ref{xy} (with $x=m-2$ and $y=1$), there is 
a $(G_2, \Sigma,3,1)$-DF, say $\cT'$, where $\Sigma$ is a $\{2^{2^{m-2}-1},3\}$-PS of $G_2$. Therefore, $\cT\,\cup\,\cT'$ is a $(G, \{G_1\}\cup\,\Sigma, 3,1)$-DF.
Note that all $2$-groups of $\Sigma$ belong to $G_1$ which is an elementary abelian $2$-group containing $2^{m}-1$ subgroups of order $2$. Hence $\Sigma'=\sigma(G_1)\,\cup\, \Sigma$ is a $\{2^{2^m-1}, 3\}$-PS and by Remark \ref{rem}, $\cT\,\cup\,\cT'$ is a DF relative to $\Sigma'$, thus completing the proof.
\end{proof}

\begin{prop}\label{ml=even1}
    There exists a $(\Z^{m-1}_{2} \times \Z_{4}\times \Z_d,\{2^{2^m-1}, 3\},3,1)$-DF whenever $m\geq 4$ is even and $d\equiv 3 \pmod{6}$, with $d\neq 9$. 
    Furthermore, there is a $(\Z^{m-1}_{2} \times \Z_{4} \times \Z_3^2,\{2^{2^m-1}, 3\},3,1)$-DF.
\end{prop}
\begin{proof}
 Set $G=\Z^{m-1}_{2} \times \Z_{4}\times H$, where 
 \[H\simeq 
  \begin{cases}
    \Z_d & \text{if $9\neq d\equiv 3 \pmod{6}$},\\
    \Z_3\times\Z_3 & \text{if $d=9$}.
  \end{cases}
 \]
 By Lemmas \ref{cyclic_DF} and \ref{DM}, there exist an $(H, S, 3,1)$-DF, with $S\simeq \Z_3$, and a $(\Z^{m-1}_{2} \times \Z_{4},3,1)$-DM.  
 Then, Lemma \ref{composition} guarantees the existence of 
 a $(G, S',3,1)$-DF, say $\cT$, where $S'=\Z^{m-1}_{2} \times \Z_{4} \times S \simeq \Z^{m-1}_{2} \times \Z_{4} \times \Z_3$. 
 Furthermore, by Lemma \ref{xy} (with $x=m$ and $y=2$), there is 
 an $(S', \{2^{2^m-1},3\},3,1)$-DF, say $\cT'$. Therefore, $\cT\,\cup\, \cT'$ yields the desired DF. 
\end{proof}

\begin{prop}\label{ml=even>1} 
  There exists a $(\Z^{m-2}_{2} \times \Z_4 \times \Z_{2^{\ell}d},\{2^{2^m-1}, 3\},3,1)$-DF 
  whenever $m\geq 4$ is even, $\ell\geq 2$, and $d\equiv 3 \pmod{6}$.
\end{prop}
\begin{proof}
Set $G=\Z^{m-2}_{2} \times H$, where $H=\Z_4 \times \Z_{2^{\ell}d}$. 
Also, letting $h_1=2^{\ell-1}d$ and $h_2=2^{\ell}d/3$, we have that
\[
  H_1 :=2\Z_4\times h_1\Z_{2^{\ell}d}\simeq \Z^{2}_{2},\;\;\;\text{and}\;\;\; 
  H_2 :=\{0\}\times h_2\Z_{2^{\ell}d}\simeq Z_3.
\]
Lemma \ref{3pyr_DF} guarantees the existence of
an $(H,\{2^3,3\},3,1)$-DF or, equivalently,
an $(H,\{H_1, H_2\},3,1)$-DF.
Also, by Lemma \ref{DM}, there is a $(\Z_2^{m-2}, 3, 1)$-DM.
Then, by Lemma \ref{composition}, we get a $(G, \{G_1, G_2\}, 3,1)$-DF, 
where $G_1 = \Z_2^{m-2}\times H_1\simeq \Z_2^{m}$
and $G_2 = \Z_2^{m-2}\times H_2 \simeq \Z^{m-2}_{2} \times \Z_3$.

As in the proof of Proposition \ref{ml=even0}, the assertion follows by 
Lemma \ref{xy} (with $x=m-2$ and $y=1$) which guarantees the existence of
a $(\Z^{m-2}_{2} \times \Z_3, \Sigma ,3,1)$-DF, where $\Sigma$ is a $\{2^{2^{m-2}-1}, 3\}$-PS.
\end{proof}


\begin{lemma}\label{mld=309} 
There exists a $(\Z_2^{3}\times \Z^2_3, \{2^7, 3^2\}, 3,1)$-DF.
\end{lemma}
\begin{proof} We provide below the 10 triples of the desired DF.
\begin{align*}
T_1 &= \{(0,0,0,0,0), (0,0,1,1,0), (0,0,0,2,2)\},\\
T_2 &= \{(0,0,0,0,0), (0,0,1,1,1), (0,0,0,1,2)\},\\
T_3 &= \{(0,0,0,0,0), (0,1,1,2,1), (1,0,0,0,2)\},\\
T_4 &= \{(0,0,0,0,0), (0,1,1,0,1), (1,0,0,2,1)\},\\
T_5 &= \{(0,0,0,0,0), (1,1,1,2,1), (0,1,0,0,2)\},\\
T_6 &= \{(0,0,0,0,0), (1,1,1,0,1), (0,1,0,2,1)\},\\
T_7 &= \{(0,0,0,0,0), (1,0,1,2,1), (1,1,0,0,2)\},\\
T_8 &= \{(0,0,0,0,0), (1,0,1,0,1), (1,1,0,2,1)\},\\
T_9 &= \{(0,0,0,0,0), (0,1,0,1,0), (1,0,0,2,0)\},\\
T_{10} &= \{(0,0,0,0,0), (0,1,0,1,1), (1,0,0,2,2)\}.
\end{align*}
\end{proof}

\begin{lemma}\label{mld=319}
   There exists a $(\Z^{2}_{2} \times \Z_4 \times \Z^2_{3},\{2^{7}, 3^2\},3,1)$-DF.
\end{lemma} 
\begin{proof} 
  Let $\Z_2^2 = \{0, \alpha,\beta,\gamma\}$, where the nonzero elements are the three involutions of $\Z_2^2$; hence, $\alpha +\beta =\gamma$. Also, set
  $G=\Z_2^2 \times H$, where $H=\Z_3 \times \Z_{12}$. Notice that $G\simeq \Z^{2}_{2} \times \Z_4 \times \Z^2_{3}$.

Consider the subset $A=\{a_i, a'_i\mid 1\leq i\leq 5\}$ of $\Z_{12}$ defined as follows:
$a_i=6-i$, $a'_i=5+i$ for $i\leq 3$, otherwise $a'_i=6+i$,
and note that 
\begin{equation}
\overline{A}:=\pm\{a_i-a'_i\mid 1\leq i\leq 5\} = \Z_{12}\setminus\{0,6\}.
\end{equation}
Now, let $\cT=\{T_1,\ldots, T_{12}\}$ be the set of triples of $G$ defined below:
  \begin{align*}
    T_i    &=\{(0,0,0), (\alpha,1,a_i), (\gamma, 2, a_i-a'_i)\},\; \text{for $1 \leq i \leq 5$,}\\
    T_{i+5}&=\{(0,0,0), (\alpha,1,a'_i), (\gamma, 2, a'_i-a_i)\},\, \text{for $1 \leq i \leq 5$,}\\
    T_{11} &=\{(0,0,0), (\alpha,1,0), (\gamma,2,0)\},\\
    T_{12} &=\{(0,0,0), (\alpha,0,4), (\gamma,0,8)\}\},\\
    T_{13} &=\{(0,0,0), (\alpha,0,1), (\gamma,0,3)\},\\
    T_{14} &=\{(0,0,0), (\alpha,0,2), (\gamma,0,11)\},\\
    T_{15} &=\{(0,0,0), (\alpha,0,9), (\gamma,0,10)\},\\        
    T_{16} &=\{(0,0,0), (0,1,4), (\alpha,1,9)\},\\
    T_{17} &=\{(0,0,0), (0,1,8), (\beta,1,3)\},\\
    T_{18} &=\{(0,0,0), (0,1,1), (\gamma,1,6)\}.              
  \end{align*}
One can check that
$\bigcup_{i=1}^{10}\Delta T_i = \bigcup_{x\in\{\alpha,\beta,\gamma\}} \{x\}\times D_x$ where
\begin{align*}
  D_\alpha = &\pm (\{1\}\times A) = \pm\big(\{1\}\times (\Z_{12}\setminus\{0,9\})\big),\\
  D_\beta = &\pm (\{1\}\times (-A)) = \pm\big(\{1\}\times (\Z_{12}\setminus\{0,3\})\big),\\
  D_\gamma = &\pm (\{1\}\times \overline{A}) = \pm\big(\{1\}\times (\Z_{12}\setminus\{0,6\})\big).
\end{align*}
Furthermore, 
$\bigcup_{i=11}^{16}\Delta T_i = \bigcup_{x\in\Z_2^2} \{x\}\times D'_x$ where
\begin{align*}
  D_0' =&\pm \{(1,1), (1,4), (1,8)\},\\
  D'_\alpha = &(\{0\}\times \Z_{12}\setminus\{0,6\})\,\cup\,\pm \{(1,0), (1,9)\},\\
  D'_\beta = &(\{0\}\times \Z_{12}\setminus\{0,6\})\,\cup\,\pm \{(1,0), (1,3)\},\\
  D'_\gamma = &(\{0\}\times \Z_{12}\setminus\{0,6\})\,\cup\,\pm \{(1,0), (1,6)\}.
\end{align*}
Therefore, 
\[\Delta \cT= 
(\{0\}\times D'_0) \,\cup\, 
(\{\alpha,\beta,\gamma\}\times \Z_3\times\Z_{12})
\setminus (\Z^2_{2} \times \{(0,0), (0,6)\}).
\]
It is left to construct a set $\cT'$ of triples such that
\begin{equation}\label{eq2:mld=319}
\begin{aligned}
\Delta \cT' =\, &(\{0\}\times \Z_3\times \Z_{12})\setminus S,\\
\text{where}\;
S=\, &(\{0\}\times D'_0) \,\cup\, \{(0,0,0), \pm (0,0,4), (0,0,6), \pm(0,1,0)\}.
\end{aligned}
\end{equation}
Indeed, $\cT\,\cup\,\cT'$ would be the desired difference family.

One can check that the set $\cT'=\{T'_1, \ldots, T'_4\}$  of the triples defined below satisfies \eqref{eq2:mld=319}.
\begin{align*}
    T'_{1} &=\{(0,0,0), (0,0,1), (0,1,11)\},\\
    T'_{2} &=\{(0,0,0), (0,0,2), (0,1,5)\},\\
    T'_{3} &=\{(0,0,0), (0,0,3), (0,1,9)\},\\
    T'_{4} &=\{(0,0,0), (0,0,5), (0,1,7)\}.         
\end{align*}
\end{proof}

\begin{lemma}\label{xy2}
Let $x\geq 3$ be an odd integer and $y=1,2$. 
Then,
  there exists a $(\Z^{x-1}_{2} \times \Z_{2^y}\times \Z^2_{3}, \{2^{2^x-1},3^2\},3,1)$-DF.
\end{lemma}
\begin{proof}
  We proceed by induction on the odd values of $x\geq 3$. If $x=3$, the result follows from Lemma \ref{mld=309} (when $y=1$) and Lemma \ref{mld=319} (when $y=2$). Assume the assertion holds for $x\geq 3$ odd, that is,
there is a $(\Z^{x-1}_{2} \times \Z_{3}\times \Z_{12}, S',3,1)$-DF, 
with $S'=H'\,\cup\, K'_1 \,\cup\, K'_2$, where 
\begin{align*}
  H' &=\Z^{x-1}_{2}\times \{0\} \times 6\Z_{12}\simeq \Z^{x}_{2},\\
  K_1' &=\{(0,0)\}\times 4\Z_{12}\simeq \Z_3,\; \text{and}\;\\
  K_2'' &=\{0\}\times \Z_{3}\times \{0\}\simeq \Z_3. 
\end{align*} 
Since, by Lemma \ref{DM}, there is a $(\Z^{2}_{2},3,1)$-DM, we use 
Lemma \ref{composition} to obtain a
$(\Z^{x+1}_{2} \times \Z_{3}\times \Z_{12}, S,3,1)$-DF, say $\cU$,
with $S=H\,\cup\, K_1 \,\cup\, K_2$, where 
\begin{align*} 
   H &= \Z_{2}^2\times H'\simeq \Z^{x+2}_{2},\\
   K_1 &= \Z_{2}^2 \times K_1'\simeq \Z_{2}^2 \times \Z_3,\; \text{and}\;\\
   K_2 &= \Z_{2}^2 \times K_2'\simeq \Z_{2}^2 \times \Z_3. 
\end{align*}
Finally, by Lemma \ref{3pyr_DF}, there is a 
$(K_i,(\Z_2^2\times\{0\})\cup(\{0\}\times K'_i),3,1)$-DF, 
say $\cU_i'$, for $i=1,2$.
Therefore, $\cU\,\cup\,\cU_1' \,\cup\,\cU_2'$ is a 
\[(\Z^{x+1}_{2} \times \Z_{3}\times \Z_{12}, \{2^f,3^2\},3,1)\text{-DF},\] with $f=2^{x+2}-1$, thus proving the induction step.
\end{proof}

\begin{prop}\label{ml=odd0}
  There exists a $(\Z^{m}_{2} \times \Z_{3} \times \Z_{d/3}, \{2^{2^m-1},3^2\},3,1)$-DF whenever $m\geq 3$ is odd and $d\equiv 9 \pmod{18}$, with $d\neq 27$. 
  Furthermore, there is a $(\Z^{m}_{2} \times \Z^3_{3},\{2^{2^m-1}, 3^2\},3,1)$-DF.
\end{prop}
\begin{proof}
Let $d=9u$, and set 
$G=\Z^{m}_{2} \times \Z_3 \times  H$, where
 \[H\simeq 
  \begin{cases}
    \Z_{3u} & \text{if $u\neq 3$},\\
    \Z_3\times\Z_3 & \text{if $u=3$}.
  \end{cases}
 \]
By Lemma \ref{cyclic_DF}, there exists a
$(H, S,3,1)$-DF where $S\simeq \Z_3$. 
By Lemma \ref{DM}, there exists a $(\Z_2^{m}\times \Z_3, 3, 1)$-DM.
Then, Lemma \ref{composition} guarantees the existence of a 
$(G, S', 3,1)$-DF, say $\cT$, where 
\[S'=\Z_2^{m}\times \Z_3\times S \simeq \Z_2^{m}\times \Z^2_3.\] 
Letting $\cT'$ be an $(S',\{2^{2^m-1}, 3^2\}, 3,1)$-DF, whose existence is guaranteed by Lemma \ref{xy2} (with $x=m$ and $y=1$), it follows that $\cT\,\cup\, \cT'$ is the desired DF.
\end{proof}

\begin{prop}\label{ml=odd1}
    There exists a $(\Z^{m-1}_{2} \times \Z_{12} \times \Z_{d/3},\{2^{2^m-1}, 3^2\},3,1)$-DF for every odd $m\geq 3$ and every $d\equiv 9 \pmod{18}$, with $d\neq 27$. 
  Furthermore, there is a 
  $(\Z^{m-1}_{2} \times \Z_{12} \times \Z^2_3, \{2^{2^m-1}, 3^2\},3,1)$-DF.
\end{prop}
\begin{proof}
Let $d=9u$, with $u$ odd, and set $G=\Z^{m-1}_{2} \times \Z_{12} \times H$, where
 \[H\simeq 
  \begin{cases}
    \Z_{3u} & \text{if $u\neq 3$},\\
    \Z_3\times\Z_3 & \text{if $u=3$}.
  \end{cases}
 \]
By Lemma \ref{cyclic_DF}, there exists a
$(H, S,3,1)$-DF where $S\simeq \Z_3$. 
By Lemma \ref{DM}, there exists a $(\Z^{m-1}_{2} \times \Z_{12},3,1)$-DM.
Then, Lemma \ref{composition}, guarantees the existence of 
 a $(G, S',3,1)$-DF, say $\cT$, where 
 \[S'=\Z^{m-1}_{2} \times \Z_{12} \times S \simeq 
 \Z^{m-1}_{2} \times \Z_{12} \times \Z_{3} \simeq  \Z^{m-1}_{2} \times \Z_4 \times \Z_3^2.\]
It is then left to build an $(S', \{2^{2^m-1},3^2\},3,1)$-DF, say $\cT'$; indeed, $\cT\,\cup\, \cT'$ would yield the desired difference family. The existence of $\cT'$ is guaranteed by 
Lemma \ref{xy2} (with $x=m$ and $y=2$) and this completes the proof.
\end{proof}

\begin{prop}\label{ml=odd>1}
 There exists a $(\Z^{m-1}_{2} \times \Z_3 \times \Z_{2^{\ell+1}d/3},\{2^{2^m-1}, 3^2\},3,1)$-DF whenever $\ell\geq 2$, $m\geq 3$ is odd and $d\equiv 9 \pmod{18}$.
\end{prop}
\begin{proof}
Set 
$G=\Z^{m-1}_{2} \times H$, where $H=\Z_3 \times \Z_{2^{\ell+1}d/3}$. By Lemma \ref{1rot_DF_abelian}, there exists an
$(H,\Sigma,3,1)$-DF where $\Sigma=\{H_1, K_1, K_2\}$ is a partial spread of $H$,
 with $|H_1|=2$ and $|K_1|=|K_2|=3$, hence of type $\{2^{1}, 3^2\}$. By Lemma \ref{DM}, there exists a $(\Z_2^{m-1}, 3, 1)$-DM.
Then, by Lemma \ref{composition}, we get a $(G, \Sigma', 3,1)$-DF, 
where $\Sigma'=\{\Z^{m}_{2}\}\,\cup\, \{\Z^{m-1}_{2} \times K_i\mid i=1,2\}$.
By Lemma \ref{xy} (with $x=m-1$ and $y=1$), there is a 
$(\Z^{m-1}_{2} \times K_i, \{2^{2^{m-1}-1}, 3\},3,1)$-DF for every $i=1,2$, and this completes the proof.
\end{proof}

\section{Proof of Theorem \ref{main}}
\label{section:main}
For the reader's convenience, we recall the statement of the main result of this paper.\\

\noindent
\textbf{Theorem \ref{main}.} 
\emph{There exists an $f$-pyramidal STS$(v)$, with $3<f<v$, over some abelian group if and only if $f=2^m-1$ for some $m\geq 3$ and one of the following conditions hold:
  \begin{enumerate}[label=$(\alph*)$]
    \item \label{main_1} $v \equiv 2^{m+1} -1 \pmod{2^m3}$,
    \item \label{main_2} $m$ is even and
    $v \equiv 2^{m} -1 \pmod{2^m3}$,
    \item \label{main_3} $m$ is odd and
    $v \equiv 2^{m} -1 \pmod{2^m9}$.
  \end{enumerate}
}

\begin{proof} Let $f>3$ and assume there is an $f$-pyramidal STS$(v)$ over an abelian group $G$. By Theorem  \ref{DF}, $G$ has order $v-f$ and it contains exactly $f<\frac{v}{2}$ involutions. Furthermore, 
there exists a $(G,\Sigma,3,1)$-DF, say $\cF$, where $\Sigma$ is a $\{2^f,3^e\}$-PS
for some $e\geq 0$; by recalling that $|\Delta \cF| = 6|\cF|$, it follows that 
\[ 
\begin{aligned}
 |\Delta \cF| &= \textstyle{\left|G\setminus \bigcup_{S\in \Sigma} S\right| = |G| - \left|\bigcup_{S\in \Sigma} S\right|} \\ 
 &= (v-f) - (1+f+2e) = v -2f -2e-1 \equiv 0 \pmod{6}.
\end{aligned}
\]
By Lemma \ref{lemma_1}.(1) and (2), we have that $f=2^m-1$ for some $m\geq 3$ (since $f>3$)
and $v-2^m+1 = v-f = |G|=2^{m}u$, for some $u>1$ (since $f< \frac{v}{2}$). Hence,
$v=2^{m}u + 2^m-1$ and 
\[
\begin{aligned}
  v -2f -2e-1 &= (2^{m}u + 2^m-1) - 2(2^m-1) - 2e -1 \\
              &= 2^{m}(u-1) - 2e \equiv 0\pmod{6},
\end{aligned}
\]
that is, 
\begin{equation}\label{main:onlyifpart}
  2^{m-1}(u-1) \equiv e\pmod{3}
\end{equation}
If $u\not\equiv 0 \pmod{3}$, then $|G|=2^mu$ is not divisible by $3$, that is, $G$ has no element of order $3$, hence $e=0$. By \eqref{main:onlyifpart}, we then have that $u\equiv 1 \pmod{3}$, say $u=3n+1$, thus $v=2^{m}u + 2^m-1 = 2^m 3n + 2^{m+1} -1$, hence 
\ref{main_1} holds. \\
If $u\equiv 0 \pmod{3}$, say $u=3n$, then $v=2^{m}3n + 2^m-1$, hence 
\ref{main_2} holds. Note that in this case, 
\eqref{main:onlyifpart} implies that $2^m\equiv e \pmod{3}$. Therefore, if $u\equiv 0 \pmod{3}$ and $m$ is odd, then $e \equiv 2^m \equiv 2 \pmod{3}$, hence $e\geq 2$. It follows that $\Sigma$ contains at least 2 distinct subgroups of $G$ of order $3$. Since $G$ is abelian, this implies that $G$ has a subgroup isomorphic to $\Z_3^2$, hence $|G|=2^mu$ is divisible by $9$. Therefore, $u\equiv 0 \pmod{9}$, say $u=9n$, thus $v = 2^{m}9n + 2^m-1$,
hence  \ref{main_3} holds.\\ 

Conversely, assume that $f=2^m-1$ for some $m\geq 3$, and one of the conditions 
\ref{main_1}-\ref{main_3} holds. Since in each of these cases $v-f \equiv 0\pmod{2^m}$,
we can write $v-f=2^{m+\ell}d$, for some $\ell\geq 0$ and some odd $d\geq 1$. \\
If \ref{main_1} holds, then $v = 2^m 3n + 2^{m+1} -1$, for some $n\geq 0$, therefore
\[2^{m+\ell}d = v-f= (2^m 3n + 2^{m+1} -1) - (2^m-1) = 2^m (3n+1),\] 
hence,
$2^\ell d=3n+1$. Since $d$ is odd, it follows that 
$d \equiv (-1)^\ell \pmod{6}$. \\
If \ref{main_2} holds, then $m$ is even  and $v = 2^m 3n + 2^m -1$, for some $n\geq 1$
(note that $n\neq 0$, since $f<v$), therefore
\[2^{m+\ell}d = v-f= 2^m 3n + 2^m -1 - (2^{m} -1)= 2^m 3n,\] 
hence, $2^\ell d=3n$. Again, since $d$ is odd, it follows that 
$d \equiv 3 \pmod{6}$.\\
Finally, if \ref{main_3} holds, then $m$ is odd and $v = 2^m 9n + 2^m -1$, for some $n\geq 1$. Therefore,
$2^{m+\ell}d = v-f= 2^m 9n +2^m -1 - (2^{m} -1)= 2^m 9n$,
hence, $2^\ell d=9n$, that is, $d \equiv 9 \pmod{18}$.

Then, Theorem \ref{main:2} guarantees the existence of a $(G,\{2^f,3^e\},3,1)$-DF, for some $e\geq 0$, where $G$ is an abelian group of order $2^{m+\ell}d = v-f$ having exactly $f=2^{m}-1$ involutions. Since $f\equiv 1,3 \pmod{6}$, by Theorem \ref{DF}, we obtain the existence of an $f$-pyramidal $STS(v)$ under the action of $G$, and this completes the proof.
\end{proof}

\section*{Acknowledgments}
Much of this research was undertaken during a visit by T. Traetta to Beijing Jiaotong University.
He expressed his sincere thanks to the 111 Project of China (B16002) for financial support and to the School of Mathematics and Statistics  at Beijing Jiaotong University for their kind hospitality. The authors' research received support from the following sources. Y.~ Chang: NSFC grant 12371326; T.~ Traetta: INDAM - GNSAGA; J.~Zhou:  NSFC grant 12171028.

We also thank the anonymous referees for their valuable comments and suggestions.

\end{document}